\documentstyle[11pt]{article}
\begin{document}
\setlength{\baselineskip}{15pt}
\title{Conservative quasipolynomial maps}
\author{Benito Hern\'{a}ndez--Bermejo $^{(1,*)}$ \and L\'{e}on Brenig $^{(2)}$}
\date{}

\maketitle

\begin{flushleft}
\noindent{\em $^{(1)}$ Escuela Superior de Ciencias Experimentales y Tecnolog\'{\i}a. Edificio Departamental II. } \\
\noindent{\em \hspace{0.5 cm} Universidad Rey Juan Carlos. Calle Tulip\'{a}n S/N. 
28933--M\'{o}stoles--Madrid. Spain.} \\
\mbox{} \\
\noindent{\em $^{(2)}$ Service de Physique Th\'{e}orique et Math\'{e}matique. Universit\'{e} Libre de Bruxelles.} \\
\noindent{\em \hspace{0.5 cm} Campus Plaine -- CP 231. Boulevard du Triomphe. B-1050 Bruxelles. Belgium.}
\end{flushleft}

\mbox{}

\mbox{}

\mbox{}

\begin{center} 
{\bf Abstract}
\end{center}
\noindent
The existence of conservative quasipolynomial (QP) maps is investigated. A classification is given for dimensions two and three, and the analytical solution of the former case is constructed. General properties of $n$-dimensional QP conservative maps are also analyzed, including their dimensional reduction and the characterization of conserved quantities.

\mbox{}

\mbox{}

\mbox{}

\mbox{}


% 03.20.+i = Classical Mechanics of Discrete Systems: general mathematical aspects.
% 03.65.Fd = Algebraic methods.
% 46.10.+z = Mechanics of discrete systems.

\mbox{}


\mbox{}

\mbox{}

\mbox{}

\mbox{}

\mbox{}

\mbox{}

\mbox{}

\mbox{}

\mbox{}

\mbox{}

\mbox{}

\noindent $^{(*)}$ {\bf Corresponding author.} 
Telephone: (+34) 91 488 73 91. Fax: (+34) 91 488 73 38. \newline 
\mbox{} \hspace{4.8cm} E-mail: {\tt bhernandez@escet.urjc.es }

\pagebreak
\begin{flushleft}
{\bf 1. Introduction}
\end{flushleft}

In the last years the study of Lotka-Volterra (LV from now on) discrete-time systems has deserved a significant attention in the literature. As in the case of their continuous (i.e. differential) counterpart \cite{lot1,vol1}, such interest has been initiated in the domain of biomathematical modelling \cite{may1}-\cite{lst1} and population dynamics 
\cite{kr1}-\cite{lw1}. Subsequently, the applications of LV maps have extended into other fields related to nonlinear dynamics such as physics \cite{raj1}-\cite{ur2}, chemistry \cite{gb1} or economy \cite{doh1}. 

In a previous article \cite{bl1} a new family of maps termed quasipolynomial (QP in what follows) was introduced. In such work it was noted that the interest of QP maps is twofold. 
First, they constitute a wide generalization of the well-known LV maps. In fact, LV maps are not just a particular QP case, but play a central (actually canonical) role in the theory of QP maps \cite{bl1}. Secondly, they are a mathematically natural discrete-time analog of the continuous QP systems, which have been extensively used in many different mathematical and applied contributions (see \cite{bl1} for a bibliography on continuous QP systems and their applications as well as for a detailed analysis of the connection between the discrete and the continuous QP formalisms). 

A classification of QP maps (or even of Lotka-Volterra maps) is still an open issue. However, some preliminary results were already demonstrated in \cite{bl1,bl2}. In this sense, a physically important kind of discrete-time systems is the conservative one (for instance, see \cite{jac} for a review). Conservative maps may arise as system models or as 
Poincar\'{e} return maps of continuous-time dynamical systems. Their importance is due to the fact that the area ($n$-volume in general) conservation property of such maps can be related to important dynamical features of differential flows such as the Liouville theorem. This explains the close relationship between conservative maps and Hamiltonian flows \cite{tab}. It is therefore not surprising that conservative maps are of central importance for the study of very diverse problems that range from the many-body problem to plasma physics, and that paradigmatic discrete-time systems such as the H\'{e}non map \cite{hen1} or the standard map \cite{chi1} belong to this class. Additionally, the property of area ($n$-volume) conservation has the relevant dynamical consequence that conservative maps cannot have attracting sets of points. Actually, this does not prevent that such maps display complex dynamics (such as chaos) compatible with the area (or volume) preservation. From an operational point of view, it is worth recalling that a map is conservative if and only if the determinant of its Jacobian matrix has value $1$ or $(-1)$ everywhere in the domain of interest. Moreover, it can be seen also that the map is orientation preserving if and only if the value of such determinant is $1$. 

The purpose of this work is to consider the existence, classification and properties of conservative QP maps. Actually, in what follows a classification of QP conservative maps will be given in the cases of dimensions 2 and 3. Moreover, in the two-dimensional case it will be possible to make use of the QP formalism in order to demonstrate that all QP conservative maps have an analytical solution which will be explicitly constructed. The $n$-dimensional situation will be also investigated, and in such case it will be demonstrated that the QP property always allows the reduction of such maps to dimension $(n-1)$.

The structure of the article is the following. In Section 2 some basic properties regarding the QP formalism for maps are reviewed in order to make the article self-contained. Section 3 is devoted to the complete classification of QP conservative maps of dimension 2 as well as to the construction of their analytical solution. In Section 4 the classification of the three-dimensional case is investigated and use of the QP formalism is made in order to reduce such maps to two-dimensional ones. In Section 5 some general properties of $n$-dimensional QP conservative maps are characterized and a general reduction procedure to dimension 
$(n-1)$ is developed. Finally, some concluding remarks are discussed in Section 6. 

\mbox{}

\mbox{}

\begin{flushleft}
{\bf 2. Overview of the QP formalism for maps}
\end{flushleft}

The aim of this section is to present an overview of the QP formalism for maps. The reader is referred to \cite{bl1} for the full details. QP maps are those of the form
\begin{equation}
\label{qpm}
	x_i(t+1)=x_i(t) \exp \left( \lambda _i + \sum _{j=1}^m A_{ij} \prod _{k=1}^n 
	[x_k(t)]^{B_{jk}} \right)
	\:\: , \:\:\:\:\: i =1, \ldots , n
\end{equation}
where {\em (i)} $m$ is an integer not necessarily equal to $n$; {\em (ii)} index $t$ is an 
integer denoting the discrete time; {\em (iii)} variables $x_i(t)$ are assumed to be positive for $i=1, \ldots ,n$ and for every $t$; and {\em (iv)} $A=(A_{ij})$, $B=(B_{ij})$ and 
$\lambda = ( \lambda _i)$ are real matrices of dimensions $n \times m$, $m \times n$ and $n \times 1$, respectively. Note that it is implicit in this definition that matrix $A$ should not have a column of zeros, and that matrix $B$ should not have a row of zeros either. The terms 
\[
	Q_j(x)=\prod _{k=1}^n [x_k(t)]^{B_{jk}} \: , \:\:\: j=1, \ldots ,m
\]
appearing in the exponential of equation (\ref{qpm}) are known as quasimonomials, and quasipolynomials are defined as linear combinations of quasimonomials. It is also convenient to introduce an additional matrix, denoted by $M$, which is of dimension $n \times (m+1)$ and is defined as:
\[
	M \equiv ( \lambda \mid A) =
	\left( \begin{array}{cccc}
		\lambda _1 & A_{11} & \ldots  & A_{1m} \\
		\vdots     & \vdots & \mbox{} & \vdots \\
		\lambda _n & A_{n1} & \ldots  & A_{nm} 
	\end{array} \right)
\]

Notice that LV maps 
\begin{equation}
\label{lvmac}
	x_i(t+1)=x_i(t) \exp \left( \lambda _i + \sum _{j=1}^n A_{ij} x_j(t) \right)
	\:\: , \:\:\:\:\: i =1, \ldots , n
\end{equation}
are a particular case of QP map, namely the one in which $B$ is the $n \times n$ identity matrix.

An important basic property is that the positive orthant is an invariant set for every QP map. This is natural in many domains (such as population dynamics) in which the 
system variables are positive by definition. In the QP context, this feature is always present. Thus, in what follows it is always assumed that QP maps are defined in int$\{ I \!\! R^n_+\}$.

A key set of transformations relating QP maps are the quasimonomial transformations (QMTs from now on) defined as: 
\[
	x_i(t)= \prod_{j=1}^n [y_j(t)]^{C_{ij}} \:\: , \:\:\: i = 1, \ldots ,n \:\: ; \:\:\: 
	\mid C \mid \neq 0
\]

The form-invariance of QP maps after a QMT is one of the cornerstones of the formalism. Actually, if we consider a $n$-dimensional QP map of matrices $A$, $B$, $\lambda$ (and $M$) and perform a QMT of matrix $C$, the result is another $n$-dimensional QP map of matrices $A'$, $B'$, $\lambda '$ (and $M'$) where:
\[
	A' = C^{-1} \cdot A \: , \:\:\: 
	B' = B \cdot C \: , \:\:\: 
	\lambda ' = C^{-1} \cdot \lambda \: , \:\:\: 
	M' = C^{-1} \cdot M 
\]

Moreover, every QMT relating two QP maps is a topological conjugacy. Consequently, we not only have a formal invariance between QP systems related by a QMT, but actually a complete dynamical equivalence (in the topological sense). These properties imply that the set of all QP maps related by means of QMTs actually constitute an equivalence class. One important label of such classes is given by the matrix products $B \cdot A$ and $B \cdot \lambda$ (or briefly $B \cdot M$) which are invariant for every equivalence class.

\mbox{}

\mbox{}

\begin{flushleft}
{\bf 3. Conservative QP maps in dimension 2}
\end{flushleft}

We proceed directly to consider two-dimensional conservative QP maps, since it can be easily demonstrated (and is intuitively clear as well) that the only one-dimensional QP conservative map is the identity, namely $x(t+1)=x(t)$. The proof is left to the reader. Consequently, we start our discussion directly with the case $n=2$. The main result of this section is the following one:

\mbox{}

\noindent{\bf Theorem 1.} {\em A QP map (\ref{qpm}) with $n=2$ is conservative if and only if the following three conditions hold:
\begin{description}
\item[\mbox{\em (1)}] $\lambda_1 + \lambda_2=0$.
\item[\mbox{\em (2)}] $A_{1i} + A_{2i} = 0$ for all $i=1, \ldots ,m$.
\item[\mbox{\em (3)}] $B_{i1}=B_{i2}$ for all $i=1, \ldots ,m$.
\end{description}
}

\mbox{}

\noindent{\bf Proof.} We denote by $J$ the Jacobian matrix of the map, with $J_{ij} = \partial F_i/\partial x_j \equiv \partial_j F_i$, where
\[
	F_{i} = x_i \exp \left( \lambda _i + \sum _{j=1}^m A_{ij} 
	\prod _{k=1}^2 x_k^{B_{jk}} \right) \equiv x_i \exp{(\varphi_i)} 
	\:\:\: , \:\: i = 1,2
\]
Let $\Delta = \mid J \mid$. It can be seen that:
\[
	\Delta = ( 1 + x_1 \partial _1 \varphi_1 + x_2 \partial _2 \varphi_2 + x_1x_2 
	[ (\partial _1 \varphi_1 )( \partial _2 \varphi_2 ) - 
	(\partial _2 \varphi_1 )(\partial _1 \varphi_2)]) \exp{(\varphi_1+\varphi_2)}
\]
Now if $\Delta = \pm 1$ everywhere in int$\{I\!\!R^2_+\}$, it is evident that Condition 2 of Theorem 1 must hold. At the same time, Condition 2 implies $\partial _1 \varphi_1=-\partial _1 \varphi_2$ and $\partial _2 \varphi_2=-\partial _2 \varphi_1$. Therefore the condition $\Delta = \pm 1$ takes the form:
\[
	\left( 1 + \sum _{j=1}^m A_{1j}(B_{j1}-B_{j2})x_1^{B_{j1}}x_2^{B_{j2}} \right) 
	\exp{(\lambda_1+	\lambda_2)}= \pm 1
\]
It is then clear that Condition 3 of the theorem must be verified, and this implies also the validity of Condition 1. The result is conservativity with $\Delta = +1$. \hfill 
{\Large $\Box$}

\mbox{}

The previous result implies a number of interesting corollaries for two-dimensional conservative QP maps:

\begin{itemize}
\item They are orientation preserving (note that $\Delta = +1$ without exception).

\item They always have a quasimonomial constant of motion of the form $I(x_1,x_2)=x_1x_2$. This is a consequence of Conditions 1 and 2 of Theorem 1 that imply:
\begin{equation}
\label{cm2d}
	\ln \left( \frac{x_1(t+1)}{x_1(t)} \right) + \ln \left( \frac{x_2(t+1)}{x_2(t)} \right) 
	=0
\end{equation}
The fact that $x_1x_2$ is a constant of motion holds from relation (\ref{cm2d}). 

\item The property of being conservative is not class-invariant in general, since the specific form of matrices $M$ and $B$ determined by Theorem 1 is usually lost after a QMT.

\item There are families of QMTs such as $C= \mbox{diag}(\gamma,\gamma)$, $\gamma \neq 0$, that preserve the conservative character of QP maps in dimension 2. This implies that if an equivalence class contains one conservative map, then such class contains an infinity of conservative maps.
\end{itemize}

The previous list of consequences can be finished with a relevant one: 

\mbox{}

\noindent{\bf Corollary 1.} {\em There are no conservative LV maps of dimension 2.}

\mbox{}

The above results also allow the development of an additional application of the QP formalism in the framework of conservative maps, namely the explicit determination of their analytical solution: 

\mbox{}

\noindent{\bf Theorem 2.} {\em The solution of every conservative two-dimensional QP map is of the form
\begin{equation}
\label{qpcmsol}
	\begin{array}{rcl}
		x_1(t) & = & x_1(0)k^t \\
		x_2(t) & = & x_2(0)k^{-t}
	\end{array}
\end{equation}
where $k$ is a positive real constant.}

\mbox{}

\noindent{\bf Proof.} Let $\{ x_1, x_2 \}$ be the variables and $A$, $B$ and $\lambda$ the matrices of the QP conservative map. We perform a QMT of matrix $C$, where
\begin{equation}
\label{cth2}
	C = C^{-1} = \left( \begin{array}{rr} 1 & 1 \\ 0 & -1 \end{array} \right)
\end{equation}
The result is a (nonconservative) QP map of variables $\{ y_1, y_2 \}$ and matrices:
\[
	M' = C^{-1} \cdot M = \left( \begin{array}{cccc} 0 & 0 & \ldots & 0 \\ 
		\lambda_1 & A_{11} & \ldots & A_{1m} \end{array} \right) \:\: , \:\:\: 
	B' = B \cdot C = \left( \begin{array}{cc} B_{11} & 0 \\ \vdots & \vdots \\ B_{m1} & 0 
		 \end{array} \right) 
\]
Then the transformed map can be written in the form:
\[
	\begin{array}{rcl}
		y_1(t+1) & = & y_1(t) \\
		y_2(t+1) & = & y_2(t) \exp \left( \lambda_1 + \sum_{j=1}^m A_{1j}(y_1(0))^{B_{j1}} 			\right) \equiv k y_2(t) \:\: , \:\:\:\: k > 0
	\end{array}
\]
And the analytic solution of this map is:
\begin{equation}
\label{scqpt}
	\begin{array}{rcl}
		y_1(t) & = & y_1(0) \\
		y_2(t) & = & y_2(0) k^t
	\end{array}
\end{equation}
Expression (\ref{qpcmsol}) is then found after applying the inverse of QMT (\ref{cth2}) to solution (\ref{scqpt}). \hfill {\Large $\Box$}

\mbox{}

Therefore the conservative maps characterized in Theorem 1 can be analytically solved in a systematic way thanks to the operational framework provided by the formalism. Actually their solution is monotonic (constant if $k=1$) and excludes dynamical behaviors such as chaos. We can now proceed to consider the three-dimensional case.

\mbox{}

\mbox{}

\begin{flushleft}
{\bf 4. Conservative QP maps in dimension 3}
\end{flushleft}

Before giving a classification of conservative QP maps in the case $n=3$, two definitions are required. The first one is the following:

\mbox{}

\noindent{\bf Definition 1.} {\em Let $B$ be the $m \times 3$ matrix of exponents of a three-dimensional QP map. Let $\bar{B}_i=(B_{i1},B_{i2},B_{i3})$ be its $i$-th row. Then in what follows $B$ will be said to be non-degenerate if the following two conditions are satisfied:
\begin{description}
\item[\mbox{\em (1)}] Provided $m \geq 3$, for every combination of three different arbitrary integers $(i,j,k)$, $1 \leq i,j,k \leq m$, the identity $\bar{B}_i=\bar{B}_j+\bar{B}_k$ is never satisfied.
\item[\mbox{\em (2)}] Provided $m \geq 4$, for every combination of four different arbitrary integers $(i,j,k,l)$, $1 \leq i,j,k,l \leq m$, the identity $\bar{B}_i+\bar{B}_j=\bar{B}_k+\bar{B}_l$ is never satisfied.
\end{description}
}

\mbox{}

The condition of being non-degenerate is not very restrictive for $B$ in practice, as it can be easily appreciated (note that it is a zero-measure condition in parameter space). In addition to Definition 1, let us recall that a real matrix is said to be non-negative if its entries are non-negative. From a practical point of view, imposing non-negativity on matrix $B$ is not very restrictive either. For instance the most commonly used QP map, the LV one, does satisfy both requirements. The relevance of this fact is not only formal: for example, in the dynamical analysis of LV population models the issues of conservativity and dissipativity must be often investigated (e.g., see \cite{hhj1}). Before presenting some examples, we proceed to state the main result of the section:

\mbox{}

\noindent{\bf Theorem 3.} {\em Consider a QP map (\ref{qpm}) with $n=3$ such that its $B$ matrix is non-negative and non-degenerate. Then such map is conservative if and only if the following four conditions hold:

\begin{description}
\item[\mbox{\em (1)}] $\lambda_1 + \lambda_2 + \lambda_3=0$.
\item[\mbox{\em (2)}]  $A_{1i} + A_{2i} + A_{3i} = 0$ for all $i=1, \ldots ,m$.
\item[\mbox{\em (3)}]  $A_{1i}B_{i1}+A_{2i}B_{i2}+A_{3i}B_{i3}=0$ for all $i=1, \ldots ,m$.
\item[\mbox{\em (4)}]  $\Omega _{12kl}+\Omega _{13kl}+\Omega _{23kl}=0$ for all $(k,l)$ such that $1 \leq k,l \leq m$ and $k<l$, where
\[
	\Omega _{ijkl} \equiv 
	\left| \begin{array}{cc} A_{ik} & A_{il} \\ A_{jk} & A_{jl} \end{array} \right| \cdot 
	\left| \begin{array}{cc} B_{ki} & B_{kj} \\ B_{li} & B_{lj} \end{array} \right|
\]
\end{description}
}

\mbox{}

\noindent{\bf Proof.} We again denote by $J$ the Jacobian matrix of the map and $\Delta = \mid J \mid$. Now we have
\[
	F_{i} = x_i \exp \left( \lambda _i + \sum _{j=1}^m A_{ij} 
	\prod _{k=1}^3 x_k^{B_{jk}} \right) \equiv x_i \exp{(\varphi_i)} 
	\:\:\: , \:\: i = 1,2,3
\]
Let $\varphi_{ij} \equiv \partial _j \varphi_i$. After some algebra it can be seen that:
\begin{equation}
\label{d3d}
	\Delta = \left\{ 1 + \sum_{i=1}^3 x_i \varphi_{ii} + \sum_{i,j=1 \: ; \: i<j}^3 x_ix_j 
	\left| \begin{array}{cc} \varphi_{ii} & \varphi_{ij} \\ \varphi_{ji} & \varphi_{jj}
	\end{array} \right| + x_1x_2x_3 \left| \begin{array}{ccc} \varphi_{11} & \varphi_{12} & 
	\varphi_{13}\\ \varphi_{21} & \varphi_{22} & \varphi_{23} \\ \varphi_{31} & \varphi_{32}
	& \varphi_{33} \end{array} \right| \right\} \exp{ \left( \sum_{i=1}^3\varphi_i \right)}
\end{equation}
The expression multiplying the exponential in (\ref{d3d}) is a quasipolynomial. This leads to  Condition 2. It can be demonstrated that $\varphi_{ij} = K_{ij}x_j^{-1}$, where $K=A \cdot Q 
\cdot B$, with $Q = \mbox{diag}(Q_1(x), \ldots , Q_m(x))$ and $Q_i(x)$ the $i$-th quasimonomial. It is also easily proved that Rank($K$) $=2$ and therefore the $3 \times 3$ determinant in (\ref{d3d}) vanishes. Consequently:
\begin{equation}
\label{d3d2}
	\Delta = \left\{ 1 + \sum_{i=1}^3 K_{ii} + \sum_{i,j=1 \: ; \: i<j}^3  \left| 
	\begin{array}{cc} K_{ii} & K_{ij} \\ K_{ji} & K_{jj} \end{array} \right| \right\} 
	\exp{ \left( \sum_{i=1}^3 \lambda_i \right)}
\end{equation}
Note that $K_{11}+K_{22}+K_{33}$ is a quasipolynomial without constant terms. Moreover, if $B$ is non-negative, the $2 \times 2$ determinants in (\ref{d3d2}) are also quasipolynomials without constant terms. This shows that if the map is conservative then $\Delta =+1$. This leads to Condition 1 of the theorem. More explicitely we find
\[
	\left| \begin{array}{cc} K_{ii} & K_{ij} \\ K_{ji} & K_{jj} \end{array} \right| = 
	\sum_{k,l=1}^m A_{ik}A_{jl}(B_{ki}B_{lj}-B_{kj}B_{li})Q_kQ_l= \sum_{k,l=1 \: ; \: k<l}^m 
	\Omega _{ijkl} Q_kQ_l
\]
where $\Omega _{ijkl}$ was defined in the statement of the theorem. Therefore: 
\begin{equation}
\label{d3dd}
	\Delta = 1 + \sum_{i=1}^m(A_{1i}B_{i1}+A_{2i}B_{i2}+A_{3i}B_{i3})Q_i+
	\sum_{k,l=1 \: ; \: k<l}^m (\Omega _{12kl}+\Omega _{13kl}+\Omega _{23kl}) Q_kQ_l
\end{equation}
Since $B$ is assumed to be non-negative and non-degenerate, then it is demonstrated that $\Delta =1$ iff Conditions 1 to 4 of the theorem hold. \hfill {\Large $\Box$}

\mbox{}

It is interesting to note that the meaning of the two conditions of Definition 1 is that degeneracies of the two possible forms (namely of the types $Q_i=Q_kQ_l$ and $Q_kQ_l=Q_{k'}Q_{l'}$, respectively) are excluded in (\ref{d3dd}). As in the two-dimensional case, some properties arise also for the conservative three-dimensional QP maps characterized in Theorem 3:

\begin{itemize}
\item They are orientation preserving (again $\Delta = +1$).

\item They always have a quasimonomial constant of motion of the type $I(x_1,x_2,x_3)=x_1x_2x_3$. This is a consequence of Conditions 1 and 2 of Theorem 3. The demonstration is analogous to that of the $n=2$ case and is thus omitted.

\item The property of being conservative is not class-invariant.
\end{itemize}

Before developing further aspects of the theory regarding the $n=3$ case, it is convenient to present some examples illustrating the previous results.

\mbox{}

\noindent{\em Example 1: Lotka-Volterra maps.}

\mbox{}

Due to its importance, our first example will be of the LV type (\ref{lvmac}) with $n=3$. Note that every LV map verifies the hypotheses of Theorem 3 regarding $B$. Conditions 1 and 2 of such theorem give some direct restrictions on the parameters of $A$ and $\lambda$. Condition 3 gives $A_{11}=A_{22}=A_{33}=0$. Finally, Condition 4 is to be applied for the three combinations of $(k,l)=\{(1,2),(1,3),(2,3)\}$ which lead to three further restrictions, namely $A_{12}A_{21}=0$, $A_{13}A_{21}=0$ and $A_{12}A_{13}=0$, respectively. There are several combinations of parameters compatible with the previous requirements. A typical one could be $A_{21}=A_{12}=0$, thus leading to:
\begin{equation}
\label{lvex1}
	\lambda = \left( \begin{array}{c} \lambda_1 \\ \lambda_2 \\ -\lambda_1-\lambda_2 
			\end{array} \right) \:\: , \:\:\:\: 
	     A  = \left( \begin{array}{c} A_{13} \\ -A_{13} \\ 0 
			\end{array} \right) \:\: , \:\:\:\: 
	     B  = \left( \begin{array}{ccc} 0 & 0 & 1 \end{array} \right) 
\end{equation}
It is simple to check for the LV map (\ref{lvex1}) that $\Delta =1$. Finally, it is also immediate to verify that the map has a first integral of the form $I=x_1x_2x_3$ as expected. Both proofs are left to the reader for the sake of brevity.

%\mbox{}
\pagebreak

\noindent{\em Example 2: Quasipolynomial maps not of Lotka-Volterra type.}

\mbox{}

As a second example, we can look for conservative maps compatible with the following form of $B$: 
\begin{equation}
\label{bex2}
	B  = \left( \begin{array}{ccc} 1 & 0 & 0 \\ 0 & 1 & 0 \\ 0 & 0 & 1 \\ 1 & 1 & 1 
	\end{array} \right) 
\end{equation}
Matrix (\ref{bex2}) complies to the requirements of Theorem 3. Conditions 1 and 2 are direct to apply. Condition 3 gives again $A_{11}=A_{22}=A_{33}=0$ plus the additional equation $A_{14}+A_{24}+A_{34}=0$. And Condition 4 now splits in 6 combinations of $(k,l)=\{(1,2),(1,3),(2,3),(1,4),(2,4),(3,4)\}$. The first three give again $A_{12}A_{21}=0$, $A_{13}A_{21}=0$ and $A_{12}A_{13}=0$ respectively, while the last three are verified identically. Then there are  several parametric combinations compatible with all these conditions. A choice analogous to the one made in Example 1 is $A_{21}=A_{12}=0$. The resulting map would be:
\begin{equation}
\label{qpex2}
	\lambda = \left( \begin{array}{c} \lambda_1 \\ \lambda_2 \\ -\lambda_1-\lambda_2 
			\end{array} \right) \:\: , \:\:\:\: 
	     A  = \left( \begin{array}{cc} A_{13} & A_{14} \\ -A_{13} & A_{24} \\ 0 & 
			- A_{14} - A_{24} \end{array} \right) \:\: , \:\:\:\: 
	     B  = \left( \begin{array}{ccc} 0 & 0 & 1 \\ 1 & 1 & 1 \end{array} \right) 
\end{equation}
Notice that matrix $B$ in (\ref{qpex2}) contains only two of the four rows initially present in the $B$ matrix (\ref{bex2}). As it was the case in Example 1, the reason is that those rows are not included because all the corresponding coefficients of such quasimonomials vanish in $A$, and therefore those quasimonomials are not present in map (\ref{qpex2}). Both the existence of the constant of motion $I=x_1x_2x_3$ and the verification of the conservative condition $\Delta = 1$ are simple to establish and thus omitted.

\mbox{}

It is to be expected that, in general, it is not possible to find an analytical solution for conservative QP maps in the case $n=3$. In spite of this, it will be demonstrated that they can be always reduced to some two-dimensional QP map (not conservative, in general). This is the result presented constructively in the next:

\mbox{}

\noindent{\bf Theorem 4.} {\em The phase space of every conservative three-dimensional QP map characterized in Theorem 3 is split in a foliation of two-dimensional manifolds on which the map is topologically conjugate to a two-dimensional QP map.}

\mbox{}

\noindent{\bf Proof.} Only Conditions 1 and 2 of Theorem 3 are needed for the decomposition. Both imply that the phase space is split in a foliation of the form $x_1x_2x_3=\mbox{constant}$. Now consider a QMT of matrix $C$ such that
\begin{equation}
\label{cpth4}
    C^{-1} = \left( \begin{array}{rrr} 1 & 0 & 0 \\ 0 & 1 & 0 \\ 1 & 1 & 1 \end{array} \right) 
\end{equation}
Let $\{ y_1, y_2, y_3 \}$ be the new variables. If we apply QMT (\ref{cpth4}) we have $C^{-1} \cdot M = M'$ or:
\[
	C^{-1} \cdot \left( \begin{array}{cccc} 
	\lambda_1 & A_{11} & \ldots  & A_{1m} \\ 
	\lambda_2 & A_{21} & \ldots & A_{2m} \\ 
	-\lambda_1-\lambda_2 & -A_{11}-A_{21} & \ldots  & -A_{1m}-A_{2m} \end{array} \right) = 
	\left( \begin{array}{cccc} 
	\lambda_1 & A_{11} & \ldots & A_{1m} \\ 
	\lambda_2 & A_{21} & \ldots & A_{2m} \\ 
	    0     &    0   & \ldots &  0   \end{array} \right) 
\]
In addition, we also have $B'=B \cdot C$ with entries $(B'_{ij})$. Then $y_3=x_1x_2x_3$ is actually a constant which can be decoupled. The result is a QP map of dimension 2 and matrices:
\[
	\tilde{M}'= \left( \begin{array}{cccc} 
	\lambda_1 & \tilde{A}_{11} & \ldots & \tilde{A}_{1m} \\ 
	\lambda_2 & \tilde{A}_{21} & \ldots & \tilde{A}_{2m} \end{array} \right) \:\: , \:\:\: 
	\tilde{B}' = \left( \begin{array}{cc} 
	B'_{11} & B'_{12} \\ \vdots & \vdots \\ B'_{m1} & B'_{m2} \end{array} \right) 
\]
where $\tilde{A}_{ij}=A_{ij}(y_3(0))^{B'_{j3}}$. This completes the proof. \hfill {\Large 
$\Box$}

\mbox{}

In general, the transformation used in the last proof is not the only possible: depending on the problem, it is often feasible to achieve similar (or even better) reductions with suitable choices of the QMT to be performed. This is illustrated in what follows with the maps characterized in Examples 1 and 2. 

\mbox{}

\noindent{\em Example 3: Lotka-Volterra maps.}

\mbox{}

We turn back to the LV system (\ref{lvex1}) obtained in Example 1. As a first possibility, we consider a QMT of matrix:
\[
	C_1 = \left( \begin{array}{rrr} -1 & 0 & 0 \\ 0 & -1 & 0 \\ 1 & 1 & 1 \end{array} \right) 
\]
Let $\{y_1,y_2,y_3\}$ be the variables of the transformed map. It can be verified that after the QMT and the decoupling of the variable $y_3$ (which is actually a constant) the result is a $2$-d non-conservative QP map of matrices:
\[
	\tilde{M}'  = \left( \begin{array}{cr} - \lambda_1 & -A_{13}y_3(0) \\ 
		- \lambda_2 & A_{13}y_3(0) \end{array} \right) \:\: , \:\:\:\: 
	\tilde{B}'  = \left( \begin{array}{cc} 1 & 1 \end{array} \right) 
\]
A better possibility is:
\begin{equation}
\label{qmt2e}
	C_2 = \left( \begin{array}{rrr} 1 & -1 & -1 \\ 0 & 1 & 0 \\ 0 & 0 & 1 \end{array} \right) 
\end{equation}
This QMT leads to the following conjugate QP map:
\begin{equation}
\label{qmtt}
	M' = \left( \begin{array}{cr} 0 & 0 \\ \lambda_2 & -A_{13} \\
		- \lambda_1 - \lambda_2 & 0 \end{array} \right) \:\: , \:\:\:\: 
	B' = \left( \begin{array}{ccc} 0 & 0 & 1 \end{array} \right) 
\end{equation}
Let $\{z_1,z_2,z_3\}$ be the variables of map (\ref{qmtt}). Then we see that $z_1(t)=z_1(0)$ is a constant. Moreover, $z_3$ is decoupled and can be solved to $z_3(t)=z_3(0) \kappa^t$, where $\kappa=\exp{(- \lambda_1 - \lambda_2)}$. Thus, after substitution we actually have a 
one-dimensional map for $z_2$:
\[
	z_2(t+1) = z_2(t) \exp{(\lambda_2-A_{13}z_3(0)\kappa^t)}
\]
In this case the result is not a QP map due to the time dependence of the exponential, but the dimensional reduction achieved is twice the one predicted by Theorem 4. Here, obviously the complete analytical solution can be obtained.

\mbox{}

\noindent{\em Example 4: Quasipolynomial maps not of Lotka-Volterra type.}

\mbox{}

As a second example, we consider the conservative QP map (\ref{qpex2}) of Example 2. We choose again the QMT of matrix (\ref{qmt2e}). The result is the QP map:
\begin{equation}
\label{qpex4}
	     M' = \left( \begin{array}{ccc} 0 & 0 & 0 \\ \lambda_2 & -A_{13} & A_{24} \\ 
		- \lambda_1 - \lambda_2 & 0 & -A_{14} - A_{24} \end{array} \right) \:\: , \:\:\:\: 
	     B' = \left( \begin{array}{ccc} 0 & 0 & 1 \\ 1 & 0 & 0 \end{array} \right) 
\end{equation}
Let $\{w_1,w_2,w_3\}$ be the variables of map (\ref{qpex4}). Again it is clear that $w_1(t)=w_1(0)$ and $w_3(t)=w_3(0) \mu^t$, with $\mu = \exp{(-\lambda_1-\lambda_2-(A_{14}+A_{24})
y_1(0))}$. Consequently, the result is also a one-dimensional (not QP) map
\[
	w_2(t+1) = w_2(t) \exp{( \tilde{\lambda}_2-A_{13}w_3(0)\mu^t)}
\]
where $\tilde{\lambda}_2=\lambda_2+A_{24}w_1(0)$. Also in this case, the analytical solution is easily obtained.

\mbox{}

\mbox{}

\begin{flushleft}
{\bf 5. Conservative QP maps in dimension n}
\end{flushleft}

From the previous sections it is clear that the task of providing a classification of conservative QP maps in dimension $n$ is not trivial. However, it is possible to demonstrate several general and relevant properties of such maps. This is the aim of the present section. The main result is the following:

\mbox{}

\noindent{\bf Theorem 5.} {\em If a QP map (\ref{qpm}) with non-negative matrix $B$ is 
conservative, then the following conditions hold:
\begin{description}
\item[\mbox{\em (1)}] $\sum_{i=1}^n \lambda_i = 0$.
\item[\mbox{\em (2)}] $\sum_{i=1}^n A_{ij} = 0$ for all $j=1, \ldots ,m$.
\end{description}
}

\mbox{}

\noindent{\bf Proof.} With analogous notation to the one in previous sections, we have:
\[
	\Delta = \left| \begin{array}{cccc} 
	1+x_1\varphi_{11} & x_1\varphi_{12} & \ldots & x_1\varphi_{1n} \\ 
	x_2\varphi_{21} & 1+x_2\varphi_{22} & \ldots & x_2\varphi_{2n} \\ 
	\vdots & \vdots & \mbox{} & \vdots \\
	x_n\varphi_{n1} & x_n\varphi_{n2} & \ldots & 1+x_n\varphi_{nn} 
	\end{array} \right| \exp \left( \sum_{i=1}^n \varphi_i \right) \equiv 
	\pi(x) \exp \left( \sum_{i=1}^n \varphi_i \right)
\]
where $\pi(x)$ is a quasipolynomial. The fact that $\Delta = \pm 1$ for conservative maps implies Condition 2 of the theorem and it can thus be written $\Delta = \pi(x) \exp \left( \sum_{i=1}^n \lambda_i \right)$. Let us now focus on $\pi (x)$. Notice that 
\begin{equation}
\label{qijnd}
	x_i\varphi_{ij} = \sum _{s=1}^m A_{is}B_{sj}x_ix_j^{-1}Q_s
\end{equation}
where $Q_s$ is the $s$-th quasimonomial. It is clear that if $B$ is non-negative, then the  quasipolynomial (\ref{qijnd}) cannot contain a constant term. Moreover, it is straightforward to verify that a product of the form
\[
	\prod_{\alpha =1}^k x_{i_{\alpha}}x_{j_{\alpha}}^{-1}Q_{s_{\alpha}} \:\: , \:\:\: 
	1 \leq k \leq n
\]
is never constant. This implies that $\pi (x)$ contains only one constant term, namely the constant 1 coming from the main diagonal, i.e. $\Delta =(1+ \sigma (x)) \exp (\sum_{i=1}^n 
\lambda _i)$, where $\sigma (x)$ is a quasipolynomial without constant terms. Therefore if the map is conservative Condition 1 must also hold. \hfill {\Large $\Box$}

\mbox{}

From the demonstration just presented some general conclusions can be obtained:

\mbox{}

\noindent{\bf Corollary 2.} {\em Every $n$-dimensional conservative QP map with non-negative matrix $B$: 
\begin{description}
\item[\mbox{\em (1)}] Is orientation preserving.
\item[\mbox{\em (2)}] Has a constant of motion of the form $I= \prod_{i=1}^n x_i$.
\item[\mbox{\em (3)}] Can lose the conservative character after a QMT, namely conservativity is not a class-invariant property.
\end{description}
}

\mbox{}

The proof of Property 2 of Corollary 2 is analogous to those considered in the particular cases of dimensions 2 and 3, and therefore is omitted for the sake of conciseness. In addition, a generalization of Theorem 4 can be established constructively:

\mbox{}

\noindent{\bf Theorem 6.} {\em The phase space of every $n$-dimensional conservative QP map with non-negative matrix $B$ is split in a foliation of $(n-1)$-dimensional manifolds on which the map is topologically conjugate to a $(n-1)$-dimensional QP map.}

\mbox{}

\noindent{\bf Proof.} The foliation is given by the manifolds $\prod_{i=1}^n x_i=\mbox{constant}$. Use of Conditions 1 and 2 of Theorem 5 is to be made in order to perform a QMT of matrix $C$, where: 
\[
	C^{-1} = \left( \begin{array}{rrrr} 1 & 0 & \ldots & 0 \\ 0 & 1 & \ldots & 0\\ 
	\vdots & \vdots & \mbox{} & \vdots \\ 1 & 1 & \ldots & 1 \end{array} \right) 
\]
Such QMT leads to a conjugate QP map of variables $\{y_1 , \ldots , y_n\}$ and matrices 
\[
	M' = \left( \begin{array}{cccc} 
	\lambda_1 & A_{11} & \ldots & A_{1m} \\ 
	\vdots & \vdots & \mbox{} & \vdots \\
	\lambda_{n-1} & A_{n-1,1} & \ldots & A_{n-1,m} \\ 
	    0     &    0   & \ldots &  0   \end{array} \right) 
\]
and $B'=B \cdot C$ (of entries to be denoted $(B'_{ij})$ in what follows). Then $y_n = 
\prod_{i=1}^n x_i$ is a constant that can be decoupled leading to a QP map of dimension 
$(n-1)$ and matrices:
\begin{equation}
\label{tmpnd}
	\tilde{M}'= \left( \begin{array}{cccc} 
	\lambda_1 & \tilde{A}_{11} & \ldots & \tilde{A}_{1m} \\ 
	\vdots    &  \vdots   & \mbox{} &  \vdots \\
	\lambda_{n-1} & \tilde{A}_{n-1,1} & \ldots & \tilde{A}_{n-1,m} 
	\end{array} \right) \:\: , \:\:\: 
	\tilde{B}' = \left( \begin{array}{ccc} 
	B'_{11} & \ldots & B'_{1,n-1} \\ \vdots & \mbox{} & \vdots \\ 
	B'_{m1} & \ldots & B'_{m,n-1} \end{array} \right) 
\end{equation}
where $\tilde{A}_{ij}=A_{ij}(y_n(0))^{B'_{jn}}$.  \hfill {\Large 
$\Box$}

\mbox{}

Note that the reduced map just obtained is not conservative in general, since matrix $\tilde{M}'$ in (\ref{tmpnd}) usually does not verify the conditions established in Theorem 5 even if $\tilde{B}'$ remains non-negative. This completes the analysis on $n$-dimensional conservative QP maps.

\mbox{}

\noindent{\em Example 5: Symplectic QP maps.}

\mbox{}

An interesting family of conservative QP maps on which the previous properties can be checked is that of symplectic QP maps \cite{bl2}. In such work it is demonstrated that a QP mapping of even dimension $n=2s$ is symplectic if and only if the following conditions hold: 
\begin{description}
\item[\mbox{\rm (a)}] $A_{ij}+A_{s+i,j}=0 \:\:$ for all $i=1, \ldots , s$, and for all $j=1, \ldots , m$.
\item[\mbox{\rm (b)}] $\lambda_{i}+\lambda_{s+i}=0 \:\:$ for all $i=1, \ldots , s$.
\item[\mbox{\rm (c)}] $A_{ip}B_{pj}=A_{ip}B_{p,s+j}=0 \:\:$ for all $i \neq j $, $\: 1 \leq i,j \leq s$, and for all $p=1, \ldots , m$.
\item[\mbox{\rm (d)}] $A_{ip}(B_{pi}-B_{p,s+i})=0 \:\:$ for all $i =1, \ldots , s$, and for all $p=1, \ldots , m$.
\end{description}
In this sense it is illustrative to observe that conditions (1) and (2) of Theorem 5 are implied by the symplectic conditions (a) and (b) above. Other general results proved in this section such as Corollary 2 or Theorem 6 are of course valid in the symplectic situation, and moreover they can be significantly improved in such case. The reader is referred to \cite{bl2} for a detailed treatment. 

\mbox{}

\mbox{}

\pagebreak
\begin{flushleft}
{\bf 6. Final remarks}
\end{flushleft}

The previous results deserve some brief considerations regarding the nature of the contributions presented. Conservative QP maps can be completely classified if $n=2$ and with great generality when $n=3$, and moreover their dimensional reduction is possible for arbitrary dimension $n$. This is to some extent remarkable as far as the presence of a common behaviour such as chaos is thus discarded in the two-dimensional case, for which the analytical solution was actually constructed. Of course, more varied dynamical behaviours can be present in QP conservative maps of higher dimensions and also when the conservative context is excluded. The non-conservative case is obviously generic in the QP framework, and it is important from the point of view of many applications as well. Such situation is mostly unexplored at present, and it seems clear that future research should also focus on the non-conservative scenario. In addition, our results lead to a better understanding of the particular but important LV case (for which statements such as Corollary 1 have been demonstrated) and put in perspective a relevant particular case of conservative QP maps, namely symplectic QP maps \cite{bl2}. 

On a different level, it can be said that the results presented provide a practical illustration of the potentialities of the QP approach. For instance, the procedures displayed in the proofs of Theorems 2, 4 and 6 (and in the corresponding examples) are just particular cases of a family of analytic algorithms of wide use in the QP formalism. Such algorithms make use of certain algebraic properties of $A$, $B$ and $\lambda$ such as the rank degeneracy, and allow important simplifications of the systems studied (simplifications of which the solution or the dimensional reduction constitute typical examples). The reader interested in a general treatment of these methods is referred to \cite{bl1} for the full details.

\pagebreak

\end{document}